\documentclass[11pt]{article}

\usepackage{graphicx}
\usepackage[
 a4paper,
 left=2.5cm,
 right=2.5cm,
 top=3cm,
 bottom=3 cm,
  ]{geometry}                         
\usepackage{color}                              
\usepackage{indentfirst}                        
\usepackage{amsmath,amssymb,amsfonts,amsthm,bm} 
\usepackage{cases}                              
\usepackage{cite}
\usepackage{mathrsfs}
\allowdisplaybreaks
\newtheorem{proposition}{Proposition}[section]
\newtheorem{theorem}{Theorem}[section]
\newtheorem{definition}{Definition}[section]

\newtheorem{example}{Example}[section]
\newtheorem{remark}{Remark}[section]
\input tcilatex

\begin{document}
\title{Weak partial $b$-metric space and Nadler's theorem}
\author{Tanzeela Kanwal\thanks{%
Govt. Degree College for Women Malakwal, Malakwal, Sargodha-40100, Pakistan. Email: tanzeelakanwal16@gmail.com},
Azhar Hussain\thanks{%
Department of Mathematics, University of Sargodha, Sargodha-40100, Pakistan. Email: hafiziqbal30@yahoo.com,}}
\date{}
\maketitle
\begin{center}
\noindent\textbf{ Abstract}
\end{center}
We study the notions of weak partial $b$-metric space and weak partial Hausdorff $b$-metric space. Moreover, we intend to generalize
Nadler's theorem in weak partial $b$-metric space by using weak partial Hausdorff $b$-metric spaces. A non-trivial example to show the validity of our result is given.

\noindent{\textbf{Mathematics Subject Classification 2010}}: 55M20, 47H10\newline

\noindent{\textbf{\ Keywords}}: Multivalued mappings, $b$-metric space, Nadler's fixed point theorem.

\section{Introduction}
To solve a particular problem in distance spaces, there have been a number of generalizations of such spaces according to the requirement and its applicability. The problem of convergence of measurable functions with respect to measure leaded to
generalized the metric space in such a way that, the set considered in metric space is replaced with the space
and consequently the function $d$ is replaced with the functional $d$. The metric space with such modifications
is named as $b$-metric space. It was first introduced as quasi metric in 1989 by Bakhtin \cite{bak}. Formally, in 1993,
Czerwik \cite{cze1, cze2} introduced the notion of $b$-metric space as a generalization of usual metric space and proved
an analogue of Banach's principle in the context of $b$-metric spaces. Later on, many authors studied the existence of fixed point for single and multivalued mapping in the context of $b$-metric space (see for example \cite{abbas, kham, bori, bota, jov, cze1, cze2} and references therein).

On the other hand, Matthews \cite{matt} introduced the notion of partial metric
space as a part of the study of denotational semantics of dataflow network.
Further, Matthews proved an analogue of Banach contraction principle in
partial metric spaces. Aydi {\it et al.} \cite{aydi} introduced the notion of partial Hausdorff metric and proved
a fixed point theorem for multi-valued mappings using the partial Hausdorff metric on partial metric space. Heckmann \cite{heck} generalized partial metric space to weak partial metric space by omitting the small self-distance
axiom. Shukla \cite{shuk} generalized the concept of partial metric to partial $b$-metric and proved some fixed point results. Beg \cite{ismat} introduced the concept of almost partial Hausdorff metric and gave a generalization
of Nadler's theorem on weak partial metric spaces using almost partial Hausdorff metric.

Motivated by above literature, we intend to generalize the concept of partial weak metric space to partial weak $b$-metric space. We define the notion of $\mathcal{H}^+$-type partial Hausdorff $b$-metric and prove Nadler's theorem on weak partial $b$-metric
spaces using $\mathcal{H}^+$-type partial Hausdorff $b$-metric. We present an example to support our existence result.

\section{Preliminaries}
Consistent with Beg \cite{ismat}, notion of weak partial metric and related concepts are as follows:
\begin{definition}\cite{ismat}
Let $X$ be a nonempty set. A function $q:X\times X\to \mathbb{R}^+$ is called weak partial metric on $X$
if for all $x,y,z\in X$, following assertions hold:
\begin{description}
  \item[(WP1)] $q(x,x)=q(x,y)~~\Leftrightarrow ~~x=y;$
  \item[(WP2)] $q(x,x)\leqq(x,y);$
  \item[(WP3)] $q(x,y)=q(y,x);$
  \item[(WP4)] $q(x,y)\leq q(x,z)+q(z,y).$
\end{description}
The pair $(X,q)$ is called weak partial metric space.
\end{definition}
\begin{example}\cite{ismat}
If $X=\{[a,b]:a,b\in \mathbb{R},a\leq b\}$, then $q([a,b],[c,d])=\max[b,d]-\min[a,c]$ is a weak partial
metric.
\end{example}
For more detail on weak partial metric space and its topology, we refer \cite{ismat}.

Let $(X,{ q })$ be a weak partial metric space and $CB^{ q }(X)$ be family of all nonempty, closed and bounded subsets
of $(X,{ q })$. For $U,V\in CB^{ q }(X)$ and $x\in X$, define $${ q }(x,U)=\inf\{{ q }(x,u),u\in U\},~~\delta_{ q }(U,V)=\sup
\{{ q }(u,V):u\in U\},$$ and $$\delta_{ q }(V,U)=\sup\{{ q }(v,U):v\in V\}.$$
Now ${ q }(x,U)=0\Rightarrow{ q }^s(x,U)=0$, where ${ q }^s(x,U)=\inf\{{ q }^s(x,u),u\in U\}$.
\begin{remark}\cite{ismat}
Let $(X,{ q })$ be a weak partial metric space and $U$ be any nonempty subset of $X$, then $$u\in\overline{U}
\Leftrightarrow { q }(u,U)={ q }(u,u).$$
\end{remark}
\begin{proposition}\label{r1}\cite{ismat}
Let $(X,{ q })$ be a weak partial metric space. For any $U,V,Y\in CB^{ q }(X)$, we have the following:
\begin{description}
  \item[(i)] $\delta_{ q }(U,U)=\sup\{{ q }(u,u):u\in U\}$;
  \item[(ii)] $\delta_{ q }(U,U)\leq\delta_{ q }(U,V)$;
  \item[(iii)] $\delta_{ q }(U,V)=0\Rightarrow U\subseteq V$;
  \item[(iv)] $\delta_{ q }(U,V)\leq \delta_{ q }(U,Y)+\delta_{ q }(Y,V)$.
\end{description}
\end{proposition}
\begin{definition}\cite{ismat}
Let $(X,{ q })$ be a weak partial metric space. For $U,V\in CB^{ q }(X)$, define $$\mathcal{H}_{{ q }}^+(U,V)=\frac12\{
\delta_{ q }(U,V)+\delta_{ q }(V,U)\}. $$ The mapping $\mathcal{H}_{{ q }}^+:CB^{ q }(X)\times CB^{ q }(X)\to [0,\infty)$, is called
$\mathcal{H}_{{ q }}^+$-type Hausdorff metric induced by ${ q }$.
\end{definition}

\begin{proposition}\cite{ismat}
Let $(X,{ q })$ be a weak partial metric space. For any $U,V,Y\in CB^{ q }(X)$, we have:
\begin{description}
  \item[(wh1)] $\mathcal{H}_{{ q }}^+(U,U)\leq\mathcal{H}_{{ q }}^+(U,V)$;
  \item[(wh2)] $\mathcal{H}_{{ q }}^+(U,V)=\mathcal{H}_{{ q }}^+(V,U)$;
  \item[(wh3)] $\mathcal{H}_{{ q }}^+(U,V)\leq \mathcal{H}_{{ q }}^+(U,Y)+\mathcal{H}_{{ q } }^+(Y,V)$.
\end{description}
\end{proposition}
\begin{definition}\cite{ismat}
Let $(X, q )$ be a weak partial metric space. A multivalued mapping $T: X\to CB^{ q }(X)$ is called $\mathcal{H}_{{ q } }^+$-contraction if
\begin{description}
  \item[$(1^{o})$] there exists $k\in(0,1)$ such that $$\mathcal{H}_{{ q } }^+(Tx\backslash \{x\},Ty\backslash \{y\})\leq k q (x,y) \text{ for every } x,y\in X,$$
  \item[$(2^{o})$] for every $x\in X$, $y$ in $Tx$ and $\epsilon>0$, there exists $z$ in $Ty$ such that $$ q (y,z)\leq \mathcal{H}_{{ q } }^+(Tx,Ty)+\epsilon.$$
\end{description}
\end{definition}
Beg \cite{ismat} gave the following variant of Nadler's fixed point theorem.
\begin{theorem}\cite{ismat}
Let $(X,q)$ be a complete weak partial metric space. Then every $\mathcal{H}_{{q}}^+$-type multivalued contraction
with lipschitz constant $k<1$ has a fixed point.
\end{theorem}
\section{Weak partial $b$-metric space and fixed point result}
We now define weak partial $b$-metric space and related concepts:
\begin{definition}\label{D1}
Let $X$ be a non empty set and $s\geq 1$ be a given real number, a function ${\sigma}:X\times %
X\to\mathbb{R}^{+}$ is called weak partial $b$-metric on $X$ if for all $x,y,z\in X$, following
conditions are satisfied:
\begin{description}
  \item[(WPB1)] ${\sigma}(x,x)={\sigma}(x,y)\Leftrightarrow x=y$;
  \item[(WPB2)] ${\sigma}(x,x)\leq{\sigma}(x,y)$;
  \item[(WPB3)] ${\sigma}(x,y)={\sigma}(y,x)$;
  \item[(WPB4)] ${\sigma}(x,y)\leq s[{\sigma}(x,z)+{\sigma}(z,y)].$
\end{description}
The pair $(X,{\sigma})$ is a weak partial $b$-metric space.
\end{definition}
\begin{example}
\begin{description}
  \item[(i)] $(\mathbb{R}^+,{\sigma})$, where ${\sigma}:\mathbb{R}^+\times \mathbb{R}^+ \to \mathbb{R}^+$ is defined as
                $${\sigma}(x,y)=|x-y|^2+1\text{  for all } x,y\in\mathbb{R}^+.$$
  \item[(ii)] $(\mathbb{R}^+,{\sigma})$, where ${\sigma}:\mathbb{R}^+\times \mathbb{R}^+ \to \mathbb{R}^+$ is defined as
                $${\sigma}(x,y)=\frac{1}{2}|x-y|^2+\max{\{a,b\}}\text{  for all } x,y\in\mathbb{R}^+.$$
\end{description}
\end{example}
\begin{definition}
A sequence $\{x_n\}$ in $(X,{\sigma})$
converges to a point $x\in X$, if and only if $${\sigma}(x,x)=\lim\limits_{n\to\infty}
{\sigma}(x,x_n).$$
\end{definition}
\begin{remark}
If ${\sigma}$ is a weak partial $b$-metric on $X$, then the function ${\sigma}^s:X\times X\to \mathbb{R}^+$ given
by ${\sigma}^s(x,y)={\sigma}(x,y)-\frac 12[{\sigma}(x,x)+{\sigma}(y,y)]$, defines a $b$-metric on $X$. Further,
a sequence $\{x_n\}$ in $(X,{\sigma}^s)$ converges to a point $x\in X$, if and only if
 \begin{equation}\label{E1}
 \lim\limits_{n,m\to\infty}{\sigma}(x_n,x_m)=\lim\limits_{n\to \infty}{\sigma}(x_n,x)={\sigma}(x,x).
\end{equation}
\end{remark}
\begin{definition}
Let $(X,{\sigma})$ be a weak partial $b$-metric space. Then
\begin{description}
  \item[(1)] A sequence $\{x_{n}\}$ in $X$ is a Cauchy sequence in $(X, \sigma)$, if it is a Cauchy sequence in metric space $(X, \sigma^{s})$.
  \item[(2)] A weak partial $b$-metric space $(X, \sigma)$ is complete if and only if the metric space $(X, \sigma^{s})$ is complete.
\end{description}
\end{definition}
Let $(X,{\sigma})$ be a weak partial $b$-metric space and $CB^{\sigma}(X)$ be family of all nonempty, closed and bounded subsets
of $(X,{\sigma})$. For $U,V\in CB^{\sigma}(X)$ and $x\in X$, define $${\sigma}(x,U)=\inf\{{\sigma}(x,u),u\in U\},~~\delta_{\sigma}(U,V)=\sup
\{{\sigma}(u,V):u\in U\},$$ and $$\delta_{\sigma}(V,U)=\sup\{{\sigma}(v,U):v\in V\}.$$
Now ${\sigma}(x,U)=0\Rightarrow{\sigma}^s(x,U)=0$, where ${\sigma}^s(x,U)=\inf\{{\sigma}^s(x,u),u\in U\}$.
\begin{remark}\label{R1}
Let $(X,{\sigma})$ be a weak partial $b$-metric space and $U$ be any nonempty subset of $X$, then $$u\in\overline{U}
\Leftrightarrow {\sigma}(u,U)={\sigma}(u,u).$$
\end{remark}
\begin{proposition}\label{P1}
Let $(X,{\sigma})$ be a weak partial $b$-metric space. For any $U,V,Y\in CB^{\sigma}(X)$, we have the following:
\begin{description}
  \item[(i)] $\delta_{\sigma}(U,U)=\sup\{{\sigma}(u,u):u\in U\}$;
  \item[(ii)] $\delta_{\sigma}(U,U)\leq\delta_{\sigma}(U,V)$;
  \item[(iii)] $\delta_{\sigma}(U,V)=0\Rightarrow U\subseteq V$;
  \item[(iv)] $\delta_{\sigma}(U,V)\leq s[\delta_{\sigma}(U,Y)+\delta_{\sigma}(Y,V)]$.
\end{description}
\begin{proof}
\begin{description}
  \item[(i)] If $U\in CB^{\sigma}(X)$, then for all $u\in U$, we have ${\sigma}(u,U)={\sigma}(u,u)$ as $\overline{U}=U$. This implies
  that $\delta_{\sigma}(U,U)=\sup\{{\sigma}(u,U):u\in U\}=\sup\{{\sigma}(u,u):u\in U\}$.
  \item[(ii)] Let $u\in U$. Since ${\sigma}(u,u)\leq{\sigma}(u,w)$ for all $w\in U$, therefore we have ${\sigma}(u,u)\leq\inf\{{\sigma}(u,v)
  :v\in V\}={\sigma}(u,V)\leq\sup\{{\sigma}(u,V):u\in U\}=\delta_{\sigma}(U,V)$.
  \item[(iii)] If $\delta_{\sigma}(U,V)=0$, then ${\sigma}(u,V)=0$ for all $u\in U$. From (i) and (ii), it follows that ${\sigma}(u,u)\leq
  \delta_{\sigma}(U,V)=0$ for all $u\in U$. Hence ${\sigma}(u,V)={\sigma}(u,u)$ for all $u\in U$. By Remark \ref{R1}, we have $u\in
  \overline{V}=V$, so $U\subseteq V$.
  \item[(iv)] Let $u\in U,v\in V$ and $y\in Y$. By (WPB4), we have $${\sigma}(u,v)\leq s[{\sigma}(u,y)+{\sigma}(y,v)].$$ Since $v\in V$ is arbitrary, therefore
      $${\sigma}(u,V)\leq s[{\sigma}(u,y)+{\sigma}(y,V)]$$
      and
      $${\sigma}(u,V)\leq s[{\sigma}(u,y)+\sup_{y\in Y}{\sigma}(y,V)],$$
      so that
      $${\sigma}(u,V)\leq s[{\sigma}(u,y)+\delta_{\sigma} (Y,V)].$$
      Since $y\in Y$ is arbitrary, therefore $${\sigma}(u,V)\leq s
  [{\sigma}(u,Y)+\delta_{\sigma} (Y,V)].$$ Since $u\in U$ is arbitrary, we have $$\delta_{\sigma}(U,V)\leq s[\delta_{\sigma}(U,Y)+\delta_{\sigma}(Y,V)].$$
\end{description}
\end{proof}
\end{proposition}
\begin{definition}\label{D2}
Let $(X,{\sigma})$ be a weak partial $b$-metric space. For $U,V\in CB^{\sigma}(X)$, the mapping $\mathcal{H}_{{\sigma}_b}^+:CB^{\sigma}(X)\times CB^{\sigma}(X)\to [0,\infty)$ define by
$$\mathcal{H}_{{\sigma}_b}^+(U,V)=\frac12\{
\delta_{\sigma}(U,V)+\delta_{\sigma}(V,U)\}$$ is called
$\mathcal{H}_{{\sigma}_{b}}^+$-type Hausdorff metric induced by ${\sigma}$.
\end{definition}
The following proposition is a consequence of Proposition \ref{P1}:
\begin{proposition}\label{P2}
Let $(X,{\sigma})$ be a weak partial $b$-metric space. For any $U,V,Y\in CB^{\sigma}(X)$, we have:
\begin{description}
  \item[(whb1)] $\mathcal{H}_{{\sigma}_b}^+(U,U)\leq\mathcal{H}_{{\sigma}_b}^+(U,V)$;
  \item[(whb2)] $\mathcal{H}_{{\sigma}_b}^+(U,V)=\mathcal{H}_{{\sigma}_b}^+(V,U)$;
  \item[(whb3)] $\mathcal{H}_{{\sigma}_b}^+(U,V)\leq s[\mathcal{H}_{{\sigma}_b}^+(U,Y)+\mathcal{H}_{{\sigma}_b}^+(Y,V)]$.
\end{description}
\begin{proof}
From (ii) of Proposition \ref{P1}, we have $$\mathcal{H}_{{\sigma}_b}^+(U,U)= \delta_{\sigma}(U,U)\leq \delta_{\sigma}(U,V)\leq \mathcal{H}_{{\sigma}_b}^+(U,V).$$
Also (whb2) obviously holds by definition. Now for (whb3), from (iv) of Proposition \ref{P1}, we have
\begin{eqnarray*}
\mathcal{H}_{{\sigma}_b}^+(U,V)&=&\frac12\{\delta_{\sigma}(U,V)+\delta_{\sigma}(V,U)\}\\
                                &\leq&\frac12\{s[\delta_{\sigma}(U,Y)+\delta_{\sigma}(Y,V)]+s[\delta_{\sigma}(V,Y)+\delta_{\sigma}(Y,U)]\}\\
                                &=& s[\frac12\{s\delta_{\sigma}(U,Y)+\delta_{\sigma}(Y,U)\}+\frac12\{\delta_{\sigma}(Y,V)+\delta_{\sigma}(V,Y)\}]\\
                                &=& s [\mathcal{H}_{{\sigma}_b}^+(U,Y)+\mathcal{H}_{{\sigma}_b}^+(Y,V)].
\end{eqnarray*}
\end{proof}
\end{proposition}
\begin{definition}\label{D3}
Let $(X,\sigma)$ be a complete weak partial $b$-metric space. A multivalued mapping $T: X\to CB^{\sigma}(X)$ is called $\mathcal{H}_{{\sigma}_b}^+$-contraction if
\begin{description}
  \item[$(1')$] for every $x,y\in X,$ there exists $k\in(0,1)$ such that $$\mathcal{H}_{{\sigma}_b}^+(Tx\backslash \{x\},Ty\backslash \{y\})\leq k\sigma(x,y);$$
  \item[$(2')$] for every $x\in X$, $y$ in $Tx$ and $\epsilon>0$, there exists $z$ in $Ty$ such that $$\sigma(y,z)\leq \mathcal{H}_{{\sigma}_b}^+(Tx,Ty)+\epsilon.$$
\end{description}
\end{definition}
\section{Fixed Point Result}
We now prove our main result:
\begin{theorem}\label{TH1}
Every $\mathcal{H}_{{\sigma}_{b}}^+$-type multivalued contraction on a complete weak partial $b$-metric space $(X,\sigma)$ has a fixed point.
\begin{proof}
Let $u_0\in X$ be arbitrary. If $u_0\in Tu_0$ then $u_0$ is the fixed point. Therefore, we assume that $u_0\not\in Tu_0$. Let $u_0\neq u_1\in Tu_0$
such that $u_1\not\in Tu_1$. From $(2')$, we have $u_2\in Tu_1$ such that $u_2\neq u_1$ and $$\sigma(u_1,u_2)\leq \mathcal{H}_{{\sigma}_
{b}}^+(Tu_0,Tu_1)+\epsilon.$$
Continuing this process we get $u_{n+1}\in Tu_n$ such that $u_{n+1}\neq u_n$ and
\begin{equation}\label{sp}
\sigma(u_n,u_{n+1})\leq \mathcal{H}_{{\sigma}_{b}}^+(Tu_{n-1},Tu_n)+\epsilon.
\end{equation}
Choosing $\epsilon=\left(\frac{1}{\sqrt{k}}-1\right)\mathcal{H}_{{\sigma}_{b}}^+(Tu_{n-1},Tu_n)$ in \eqref{sp}, we have
 $$\sigma(u_n,u_{n+1})\leq \mathcal{H}_{{\sigma}_{b}}^+(Tu_{n-1},Tu_n)+\left(\frac{1}{\sqrt{k}}-1\right)\mathcal{H}_{{\sigma}_{b}}^+(Tu_{n-1},Tu_n)=
 \frac{1}{\sqrt{k}}\mathcal{H}_{{\sigma}_{b}}^+(Tu_{n-1},Tu_n).$$
 Thus $$\sqrt{k}\sigma(u_n,u_{n+1})\leq \mathcal{H}_{{\sigma}_{b}}^+(Tu_{n-1},Tu_n)=\mathcal{H}_{{\sigma}_{b}}^+\left(Tu_{n-1}\backslash
 \{u_{n-1}\},Tu_n\backslash \{u_n\}\right).$$
 From $(1')$, we get
 $$\sqrt{k}\sigma(u_n,u_{n+1})\leq k\sigma(u_{n-1},u_n)=(\sqrt{k})^2\sigma(u_{n-1},u_n).$$
 Thus for all $n\in\mathbb{N}$,
 $$\sigma(u_n,u_{n+1})\leq \sqrt{k}\sigma(u_{n-1},u_n).$$
 Inductively, we get
 $$\sigma(u_n,u_{n+1})\leq (\sqrt{k})^n\sigma(u_{0},u_1).$$
 Using (WPB4) of weak partial $b$-metric, for some $m\in \mathbb{N}$, we have
 \begin{eqnarray*}
 \sigma(u_n,u_{n+m})&\leq& s\sigma(u_n,u_{n+1})+s^2\sigma(u_{n+1},u_{n+2})+\cdot\cdot\cdot+s^{m}\sigma(u_{n+m-1},u_{n+m})\\
                    &\leq& s(\sqrt{k})^n\sigma(u_0,u_{1})+s^2(\sqrt{k})^{n+1}\sigma(u_{0},u_{1})+\cdot\cdot\cdot+s^{m}(\sqrt{k})^{n+m-1}\sigma(u_{0},u_{1})\\
                    &=& \left(s(\sqrt{k})^{n}+s^2(\sqrt{k})^{n+1}+\cdot\cdot\cdot+s^{m}(\sqrt{k})^{n+m-1}\right)\sigma(u_0,u_{1})\\
                    &=&\frac{1}{s^{n-1}}\left(s^{n} (\sqrt{k})^{n}+s^{n+1} (\sqrt{k})^{n+1}+\cdot\cdot\cdot+s^{n+m-1} (\sqrt{k})^{n+m-1}\right)\sigma(u_0,u_{1})\\
                    &=&\frac{1}{s^{n-1}}\sum_{i={n}}^{n+m-1}s^i(\sqrt{k})^{i}. \sigma(u_0,u_{1})\\
                    &<&\frac{1}{s^{n-1}}\sum_{i={n}}^{\infty}s^i(\sqrt{k})^{i}. \sigma(u_0,u_{1}),
 \end{eqnarray*}
 this implies $$\frac{1}{s^{n-1}}\sum_{i={n}}^{\infty}s^i(\sqrt{k})^{i}. \sigma(u_0,u_{1})\to 0~~ \text{as}~~ n\to +\infty.$$
 Using \eqref{E1} and the definition of $\sigma^s$, we get
 $$\sigma^s(u_n,u_{n+m})\leq 2\sigma(u_n,u_{n+m})\to 0~~ \text{ as }~~ n\to +\infty.$$
 This implies that $\{u_n\}$ is a Cauchy sequence in $b$-metric space $(X,\sigma^s)$. Since $(X,\sigma)$ is complete, therefore
 $ (X,\sigma^s)$ is a complete $b$-metric space. Consequently, the sequence $\{u_n\}$ converges to a point (say) $u^*\in X$ w.r.t $b$-metric $\sigma^s$,
 that is, $\lim\limits_{n\to+\infty}\sigma^s(u_n,u^*)=0$. Again, from \eqref{E1} we get
 $$\sigma(u^*,u^*)=\lim\limits_{n\to+\infty}\sigma(u_n,u^*)=\lim\limits_{n\to+\infty}\sigma(u_n,u_n)=0.$$
 Thus we conclude that $\{u_n\}$ is a Cauchy sequence in $(X,\sigma)$. Since $(X,\sigma)$ is complete. Therefore, there exists $u^*\in X$ such that
 $\lim\limits_{n\to+\infty}u_n=u^*$. Now, we show that $u^*\in T$, that is, a fixed point of $T$. Suppose on contrary that $u^*\not\in Tu^*$. Since
 \begin{eqnarray*}
 \frac12[\delta_{\sigma}(Tu_n,Tu^*)+\delta_{\sigma}(Tu^*,Tu_n)]&=&\mathcal{H}_{{\sigma}_{b}}^+(Tu_{n},Tu^*)\\
                                                                &=&\mathcal{H}_{{\sigma}_{b}}^+(Tu_{n}\backslash\{u_n\},Tu^*\backslash\{u^*\})\\
                                                                &\leq& k\sigma(u_n,u^*),
 \end{eqnarray*}
 it follows that
 $$\lim_{n\to+\infty}\inf[\delta_{\sigma}(Tu_n,Tu^*)+\delta_{\sigma}(Tu^*,Tu_n)].$$
 Since
 $$\lim_{n\to+\infty}\inf\delta_{\sigma}(Tu_n,Tu^*)+\lim_{n\to+\infty}\inf\delta_{\sigma}(Tu^*,Tu_n)\leq
 \lim_{n\to+\infty}\inf[\delta_{\sigma}(Tu_n,Tu^*)+\delta_{\sigma}(Tu^*,Tu_n)],$$
 we have
 $$\lim_{n\to+\infty}\inf\delta_{\sigma}(Tu_n,Tu^*)+\lim_{n\to+\infty}\inf\delta_{\sigma}(Tu^*,Tu_n)=0.$$
 This implies that
 $$\lim_{n\to+\infty}\inf\delta_{\sigma}(Tu_n,Tu^*)=0.$$
 Since
 $$\sigma(u^*,Tu^*)\leq \delta_{\sigma}(Tu_n,Tu^*)+\sigma(u_{n+1},u^*),$$
 it follows that
 \begin{eqnarray*}
 \sigma(u^*,Tu^*)&\leq&\lim_{n\to+\infty}\inf[\delta_{\sigma}(Tu_n,Tu^*)+\sigma(u_{n+1},u^*)]\\
                 &=&\lim_{n\to+\infty}\inf\delta_{\sigma}(Tu_n,Tu^*)+\lim_{n\to+\infty}\sigma(u_{n+1},u^*).
 \end{eqnarray*}
 This implies $\sigma(u^*,Tu^*)=0$, therefore from \eqref{E1}, we obtain
 $$\sigma(u^*,u^*)=\sigma(u^*,Tu^*),$$
 which implies $u^*\in \overline{Tu^*}=Tu^*$, as $Tu^*$ is closed.
\end{proof}
\end{theorem}
\begin{example}
Consider a set $X=\{0,\frac12,1\}$ endowed with weak partial $b$-metric $\sigma: X\times X\to \mathbb{R}^+$
given by
$$\sigma(u,v)=\frac12|u-v|^2+\frac12\max\{u,v\}~~\text{ for all }~~u,v\in X.$$
Since $\sigma\left(\frac12,\frac12\right)=\frac14\neq0$ and $\sigma(1,1)=\frac12\neq0$. Also
\begin{eqnarray*}
u\in \overline{\{0\}}&\Leftrightarrow& \sigma(u,\{0\})=\sigma(u,u)\\
                    &\Leftrightarrow& \frac12u^2+\frac12u=\frac12u\Leftrightarrow u=0\\
                    &\Leftrightarrow& u\in \{0\}.
\end{eqnarray*}
Hence, $\{0\}$ is closed with respect to the weak partial $b$-metric $\sigma$.
\begin{eqnarray*}
u\in\overline{\{0,1\}}&\Leftrightarrow&\sigma(u,\{0,1\})=\sigma(u,u)\\
                      &\Leftrightarrow&\min\left\{\frac12u^2+\frac12u,\frac12|u-1|^2+\frac12\max\{u,1\}\right\}=\frac12u\\
                      &\Leftrightarrow&u\in\{0,1\}.
\end{eqnarray*}
Hence, $\{0,1\}$ is closed with respect to the weak partial $b$-metric $\sigma$.
\begin{eqnarray*}
u\in\overline{\left\{0,\frac12\right\}}&\Leftrightarrow&\sigma\left(u,\left\{0,\frac12\right\}\right)=\sigma(u,u)\\
                      &\Leftrightarrow&\min\left\{\frac12u^2+\frac12u,\frac12\left|u-\frac12\right|^2+\frac12\max\left\{u,\frac12\right\}\right\}=\frac12u\\
                      &\Leftrightarrow&u\in\left\{0,\frac12\right\}.
\end{eqnarray*}
Define a mapping $T:X\to CB^{\sigma}(X)$ by
$$T(0)=\{0\},~~T\left(\frac12\right)=\left\{0,\frac12\right\}~~\text{ and }~~T(1)=\{0,1\}.$$

We show that for all $u,v\in X$, the contractive condition $(1')$ holds for all $k\in (0,1)$. For this we consider following cases:
\begin{description}
  \item[] For $u=v=0$, we have
            $$\mathcal{H}_{{\sigma}_{b}}^+\left(T(0)\backslash \{0\},T(0)\backslash \{0\}\right)
            =\mathcal{H}_{{\sigma}_{b}}^+\left(\{0\}\backslash \{0\},\{0\}\backslash \{0\}\right)
            =\mathcal{H}_{{\sigma}_{b}}^+\left(\emptyset,\emptyset\right)=0,$$
            so $(1')$ satisfied.
  \item[] For $u=0, v=\frac{1}{2}$, we have
            $$\mathcal{H}_{{\sigma}_{b}}^+\left(T(0)\backslash \{0\},T\left(\frac12\right)\backslash \left\{\frac12\right\}\right)
            =\mathcal{H}_{{\sigma}_{b}}^+\left(\{0\}\backslash \{0\},\left\{0,\frac12\right\}\backslash \left\{\frac12\right\}\right)
            =\mathcal{H}_{{\sigma}_{b}}^+\left(\emptyset,\{0\}\right)=0,$$
            so $(1')$ satisfied.
  \item[] For $u=v=\frac{1}{2}$, we have
            \begin{eqnarray*}
            &&\mathcal{H}_{{\sigma}_{b}}^+\left(T\left(\frac12\right)\backslash \left\{\frac12\right\},T\left(\frac12\right)\backslash \left\{\frac12\right\}\right)
            =\mathcal{H}_{{\sigma}_{b}}^+\left(\left\{0,\frac12\right\}\backslash \left\{\frac12\right\},\left\{0,\frac12\right\}\backslash \left\{\frac12\right\}\right)\\
            &&=\mathcal{H}_{{\sigma}_{b}}^+\left(\{0\},\{0\}\right)=\sigma(0,0)=0,
            \end{eqnarray*}
            so $(1')$ satisfied.
  \item[] For $u=0,v=1$, we have
            $$\mathcal{H}_{{\sigma}_{b}}^+\left(T(0)\backslash \{1\},T(1)\backslash \{1\}\right)
            =\mathcal{H}_{{\sigma}_{b}}^+\left(\{0\}\backslash \{0\},\{0,1\}\backslash \{0\}\right)
            =\mathcal{H}_{{\sigma}_{b}}^+\left(\emptyset,\{0\}\right)=0,$$
            so $(1')$ satisfied.
  \item[] For $u=\frac12, v=1$, we have
            \begin{eqnarray*}&&\mathcal{H}_{{\sigma}_{b}}^+\left(T\left(\frac12\right)\backslash \left\{\frac12\right\},T(1)\backslash \{1\}\right)
            =\mathcal{H}_{{\sigma}_{b}}^+\left(\left\{0,\frac12\right\}\backslash \{0\},\{0,1\}\backslash \{1\}\right)
            =\mathcal{H}_{{\sigma}_{b}}^+\left(\{0\},\{0\}\right)\\&&=\sigma(0,0)=0,\end{eqnarray*}
            so $(1')$ satisfied.
  \item[] For $u=v=1$, we have
            $$\mathcal{H}_{{\sigma}_{b}}^+\left(T(1)\backslash \{1\},T(1)\backslash \{1\}\right)
            =\mathcal{H}_{{\sigma}_{b}}^+\left(\{0,1\}\backslash \{1\},\{0,1\}\backslash \{1\}\right)
            =\mathcal{H}_{{\sigma}_{b}}^+\left(\{0\},\{0\}\right)=\sigma(0,0)=0,$$
            so $(1')$ satisfied.
\end{description}
Further, we show that for every $u\in X, v\in Tu$ and $\epsilon>0$, there exists $w\in Tv$ such that
$$\sigma(v,w)\leq\mathcal{H}_{{\sigma}_{b}}^+\left(Tu,Tv\right)+\epsilon.$$
So,
\begin{description}
  \item[(a)] If $u=0,v\in T(0)=\{0\},\epsilon >0$, there exists $w\in Tv=\{0\}$ such that
            $$0=\sigma(v,w)\leq \mathcal{H}_{{\sigma}_{b}}^+(Tv,Tu)+\epsilon.$$
  \item[(b)] If $u=\frac12,v\in Tu=T(\frac{1}{2})=\{0,\frac{1}{2}\}$, for $v= 0,\epsilon >0$, there exists $w\in Tv=\{0\}$ such that
            $$0=\sigma(v,w)<\frac{3}{16}+\epsilon\leq \mathcal{H}_{{\sigma}_{b}}^+(Tv,Tu)+\epsilon.$$
  and for $v= \frac{1}{2}, \epsilon >0$, there exists $w\in Tv=\{0,\frac12\}$ such that
            $$\frac14=\sigma(v,w)<\frac{1}{4}+\epsilon\leq \mathcal{H}_{{\sigma}_{b}}^+(Tv,Tu)+\epsilon.$$
  \item[(d)] If $u=1,v\in Tu=T(1)=\{0,1\}$, for $v= 0,\epsilon >0$, there exists $w\in Tv=\{0\}$ such that
            $$0=\sigma(v,w)<\frac{3}{4}+\epsilon\leq \mathcal{H}_{{\sigma}_{b}}^+(Tv,Tu)+\epsilon.$$
 and for $v= 1,\epsilon >0$, there exists $w\in Tv=\{0,1\}$ such that
            $$\frac{1}{2}=\sigma(v,w)<\frac{1}{2}+\epsilon\leq \mathcal{H}_{{\sigma}_{b}}^+(Tv,Tu)+\epsilon.$$
\end{description}
Thus the condition $(2')$ is also satisfied.

Hence Theorem \ref{TH1} can be applied and we conclude that $u\in \{0,\frac{1}{2} ,1\}$ is fixed points of $T$.
\end{example}
\bibliographystyle{amsplain}

\end{document}